\documentclass[12pt]{article}
\usepackage{amsmath}
\usepackage{amssymb}
\usepackage{amscd}
\usepackage{epsfig}
\usepackage{amsthm}
\usepackage{psfig} 

\def\no{\noindent}
\newtheorem{dfn}{Definition}[section]
\newtheorem{rem}[dfn]{Remark}
\newtheorem{thm}[dfn]{Theorem}

\newtheorem{defn}[dfn]{Definition}

\newtheorem{lem}[dfn]{Lemma}

\newtheorem{prob}[dfn]{Problem}
\newtheorem{prop}[dfn]{Proposition}
\newtheorem{assumption}[dfn]{Assumption}

\newtheorem{cor}[dfn]{Corollary}

\def\CR{\curvearrowright}
\def\acts{\CR}

\def\R{{\mathbb R}}
\def\eps{\epsilon}
\def\al{\alpha}
\def\be{\beta}
\def\ga{\gamma}
\def\Z{\mathbb Z}

\def\H{\mathbb H}
\def\Ga{\Gamma}

\def\Del{\Delta}
\def\Q{\mathbb Q}
\def\ol{\overline}
\def\A{\mathbb A}
\def\si{\sigma}
\def\F{{\mathcal F}}

\def\Om{{\Omega}}
\def\Si{\Sigma}
\def\la{\lambda}
\def\t{\tilde}

\begin{document}

\title{Convex projective structures on Gromov--Thurston manifolds}
\author{Michael Kapovich}

\maketitle
\begin{abstract}
We consider Gromov--Thurston examples of negatively curved
$n$-manifolds which do not admit metrics of constant sectional
curvature. We show that for each $n\ge 4$ some of the
Gromov--Thurston manifolds admit strictly convex real--projective
structures.
\end{abstract}

\section{Introduction}

Gromov and Thurston in \cite{Gromov-Thurston(1987)} constructed,
for each $n\ge 4$, examples of compact $n$-manifolds which admit
metrics of negative curvature, with arbitrarily small pinching
constants, but do not admit metrics of constant curvature. We
review these examples in section \ref{gte}. 
The main goal of this paper is to put {\em convex projective
structures} on Gromov--Thurston examples. Suppose that $\Om\subset
\R P^n$ is an open subset and $\Ga\subset PGL(n+1,\R)$ is a
subgroup acting properly discontinuously on $\Om$. The quotient
orbifold $Q=\Om/\Ga$ has natural projective structure $c$. The
structure $c$ is said to be {\em (strictly) convex} iff $\Om$ is a
(strictly) convex proper subset of $\R P^n$. In this case we refer
to $Q$ as (strictly) convex projective orbifold.

Our main result then is:

\begin{thm}\label{main}
Gromov-Thurston examples admit strictly convex projective
structures.
\end{thm}

We refer the reader to section \ref{last} for the more precise
statement. Our theorem will be proven in section \ref{last} via
``bending" of the original hyperbolic structure on a certain
hyperbolic manifold $M$ (used to construct Gromov-Thurston
examples) in the manner similar to \cite{Gromov-Thurston(1987)},
where flat-conformal structures were constructed on certain
negatively curved manifolds.

There are two parts in this proof: (1) Producing a projective
structure, (2) proving that the structure  is convex. Then strict
convexity of the structure follows from Benoist's theorem below
(theorem \ref{benoist}), since Gromov-Thurston examples have
Gromov-hyperbolic fundamental groups.

Part (1) is dealt with by solving a certain {\em product of
matrices} problem, which is a special case of a  Lie-theoretic
problem interesting on its own right, see section \ref{pro}. The
projective manifolds $M'$ are then built by gluing convex subsets
of the hyperbolic manifolds $M$. By passing to the universal cover
we obtain a tessellation of $\t{M}'$ by convex polyhedra in
$\H^n$, each of which has infinitely many facets.

Dealing with (2) is especially interesting, since, at present,
there is only one general method for proving convexity of
projective structures, namely via Vinberg--Tits fundamental domain
theorem \cite{Vinberg(1971)}. Unfortunately, this theorem applies
only to reflection groups, which cannot be used in higher
dimensions. Our approach to proving convexity is to adapt
Vinberg's arguments in a more general context of manifolds
obtained by gluing convex cones with {\em infinitely many faces}.
In this setting, Vinberg's arguments (requiring polyhedrality of
the cones) do not directly apply and we modify them by appealing
to the {\em small cancellation theory}, see section
\ref{convexitythm}.

\medskip
The main motivation for this paper comes from the following
beautiful

\begin{thm}
\label{benoist} (Y. Benoist, \cite{benoist(2004)}) Suppose that a
convex projective orbifold $M$ is compact. Then $M$ is strictly convex iff
$\Ga=\pi_1(M)$ is Gromov-hyperbolic.
\end{thm}

Examples of convex-projective structures on compact orbifolds are
provided by the quotients of  round balls in $\R P^n$ by discrete
cocompact groups of automorphisms. The Hilbert metric on such
examples is a Riemannian metric of constant negative sectional
curvature. Thus such orbifolds are {\em hyperbolic}. By deforming
the above examples in $\R P^n$ one obtains other examples of
strictly convex projective manifolds/orbifolds.

In 2002 I was asked by Bruce Kleiner and Francois Labourie if one
can construct examples of compact strictly convex projective
manifolds which are not obtained by deforming hyperbolic examples. 
The main goal of this paper is to prove that such examples indeed
exist in all dimensions $\ge 4$. Independently, such examples were
constructed by Yves Benoist in dimension $4$ using reflection
groups, see \cite{benoist(2006)}. The paper \cite{benoist(2006)}
also produces ``exotic'' strictly convex subsets $\Om$ in $\R P^n$
for all $n\ge 3$: The metric space $(\Om, d_H)$ is Gromov-hyperbolic but is not
quasi-isometric to $\H^n$, where $d_H$ is the Hilbert metric on
$\Om$. However these examples do not appear to admit discrete
cocompact groups of automorphisms.

\medskip
{\bf Acknowledgments.}  During the work  on this paper I was
partially supported by the NSF grant DMS-04-05180. I am grateful
to John Millson for explaining to me the construction of
arithmetic examples in section \ref{gte} and to Yves Benoist and Bruce Kleiner for
 useful conversations.

\section{Gromov-Thurston examples}
\label{gte}

In this section we review Gromov-Thurston
 examples \cite{Gromov-Thurston(1987)} of compact
$n$-manifolds (more generally, orbifolds) which admit metrics of negative curvature
but do not admit metrics of constant curvature. (Note that Gromov
and Thurston \cite{Gromov-Thurston(1987)} have other examples as
well: These examples will not be discussed here.)

\medskip
Consider the quadratic form
$$
\varphi(x)= x_1^2+...+x_n^2 - \sqrt{p} x_{n+1}^2
$$
where $p$ is a (positive) prime number, $n\ge 2$. Let $\t\Ga=
Aut(\varphi)\cap GL(n+1,\Z)$; then $\t\Ga$ is a cocompact arithmetic
subgroup in $Aut(\varphi)\cong O(n,1)$.

We let $H$ denote the Lorentzian model of the hyperbolic space $\H^n$:
$$
\{x: \varphi(x)=-1, x_{n+1}>0\}.
$$
Consider the linear subspace
$$
V=\{x\in \R^{n+1}: x_1=x_2=0\}.
$$
The intersection $V\cap H$ is a totally-geodesic codimension 2
hyperbolic subspace. The stabilizer of $V$ in $\t\Ga$ acts
cocompactly on $V\cap H$ since it is isomorphic to the set of
integer points in the algebraic group
$$
Aut(x_3^2+...+x_n^2 - \sqrt{p} x_{n+1}^2).
$$

Suppose that $W\subset \R^{n+1}$ is a {\em rational} codimension 1
linear subspace containing $V$. Then  the Lorentzian (with respect
to $\varphi$) involution $\tau_W$ fixing $W$ pointwise belongs to
$GL(n+1,\Q)$. Observe that the groups $\t\Ga$ and $\tau_W \t\Ga
\tau_W$ are commensurable. Therefore, there exists a finite index
subgroup $\Ga_W\subset \t\Ga$ which is normalized by $\tau_W$. By
applying this procedure to two appropriately chosen rational
hyperplanes passing through $V$ we obtain

\begin{lem}
Given a number $m\ge 1$ there exists a subgroup $\hat\Ga\subset
Aut(\varphi)$ commensurable to $\t\Ga$, which contains a dihedral
subgroup $D_m$ fixing $V$ pointwise. The generating involutions in $D_m$ acts
as reflections.
\end{lem}

By passing to an appropriate torsion-free normal subgroup $\Ga\subset \hat\Ga$ we get  a compact
hyperbolic manifold $M=H/\Ga$. Let $\bar\Ga$ denote the subgroup of $\hat\Ga$
generated by $\Ga$ and the dihedral subgroup $D_m$.

The group $D_m$ acts on $M$ isometrically
with a fundamental domain $O$ (that can be identified with the
orbifold $M/D_m=\H^n/\bar\Ga$), which is a manifold with corners
so that the corner (possibly disconnected) corresponds to the
hyperbolic subspace $V\cap \H^n$. The dihedral angle at this
corner is $\pi/m$. By abusing notation we will keep the notation
$V$ for this codimension 2 totally-geodesic submanifold of $M$.

The boundary of $O\setminus V$ is the union of two codimension 1
totally-geodesic (possibly disconnected) submanifolds. We denote
the closures of these submanifolds $W_1, W_2$: these are
submanifolds with boundary (which is equal to $V$) in $M$. Then we
can think of the manifold $M$ as obtained by gluing $2m$ copies of
$O$.

\begin{assumption}\label{ass0}
We assume from now on that the manifold $M$ admits an isometric
 action $D_{2m}\acts M$ of a dihedral group $D_{2m}$ which
 contains $D_m$ as an index 2 subgroup.
\end{assumption}

\noindent Then there is an isometric involution $\iota: O\to O$ which
interchanges $W_1$ and $W_2$. We now construct new manifolds (without boundary) $M'$ by gluing
$2m-2$ copies  of ${O}$.

\begin{figure}[tbh]
\begin{center}
{\input{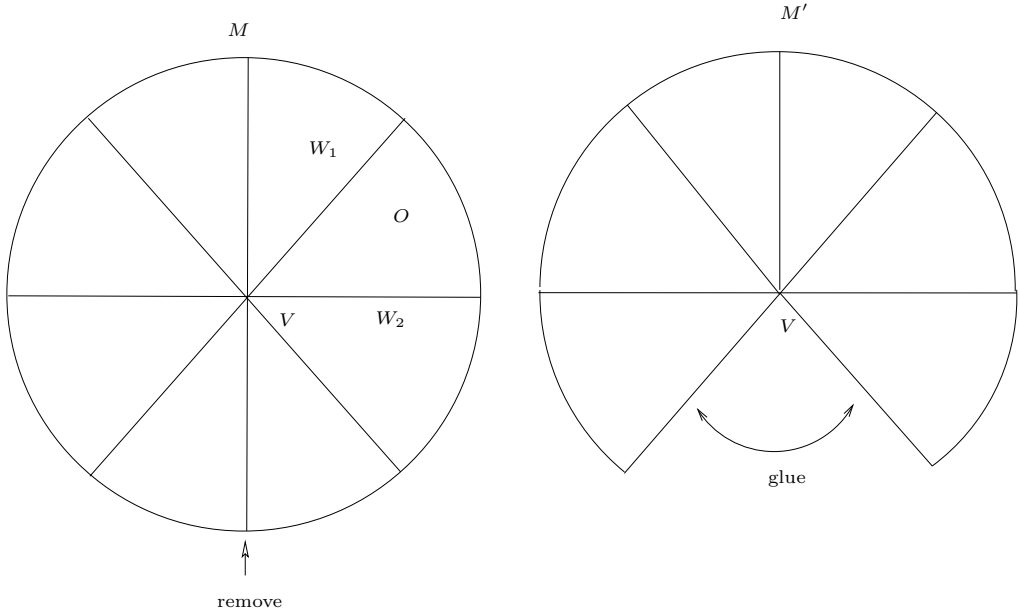}}
\end{center}
\caption{\sl Constructing the manifold $M'$. }
\label{f0.fig}
\end{figure}

\begin{rem}
Another class of Gromov--Thurston examples $M''$ is obtained by
gluing $2m+2$ copies of $O$. Construction of projective structures
on such manifolds is a bit more complicated than the one explained
in this paper, therefore the manifolds $M''$ will not be discussed
here.
\end{rem}

We will think of $M, M'$ as {\em doubles} of the manifolds $N, N'$
which are obtained by gluing $m, m-1$ copies of ${O}$
respectively. Thus $M'$  is obtained by ``subtracting'' two copies
of $O$ to $M$.

In section \ref{last} we will use an alternative description of
$M'$. Assumption \ref{ass0} implies that the manifold $N'$ admits a reflection symmetry $\theta'$ fixing
the submanifold $V$. Then $M'$   is diffeomorphic to
the manifold obtained by gluing two copies of $N'$
via the involution $\theta'|\partial N'$   of the boundary.

\begin{prop}
[Gromov, Thurston, \cite{Gromov-Thurston(1987)}]  For sufficiently
large $m$, the manifold $M'$ admits a metric of negative sectional
curvature varying in the interval $[-1+\eps_m, -1]$. Moreover,
$\lim_{m\to\infty} \eps_m=0$.
\end{prop}

\begin{rem}
Note that $M'$ admits a canonical singular Riemannian metric which
is smooth and hyperbolic away from $V$. The negatively curved
Riemannian metric on $M'$ is obtained by modifying the above
singular metric on a regular $R$-neighborhood of $V$.  This
modification works provided that $R$ is sufficiently large, which
is achieved by taking large $m$. Alternatively, one can fix $m$
and pass to an appropriate finite-index subgroup of $\Ga$.
\end{rem}

Thus $\pi_1(M')$ is Gromov-hyperbolic provided that $m$ is
sufficiently large.

\begin{prop}
[Gromov, Thurston, \cite{Gromov-Thurston(1987)}] If $n\ge 4$ then
$M'$ does not admit a metric of constant (negative) curvature.
\end{prop}
\proof Our argument is a variation on the argument given in \cite{Gromov-Thurston(1987)}.
The idea of the proof is to apply
Mostow Rigidity Theorem several times both in dimension $n$ and $n-1$.

Suppose that $M'$ admits a hyperbolic metric $g$. Observe that the
group $F=D_{m-1}$ acts (via homeomorphisms) on $M'$. Therefore, by
Mostow Rigidity Theorem, $F\acts M'$ is homotopic to an isometric
action $F\acts (M',g)$; the fixed-point set of this action is a
submanifold $V'$ homotopic to $V$.

The fundamental domain for the latter action is a submanifold with
boundary $O'\subset M'$ homotopic to $O$. Thus the dihedral angle
of $O'$ along $V'$ equals $\frac{\pi}{m-1}$.

The group $F\acts M'$ contains a topological reflection $\si$
fixing $S:=\partial N'$ pointwise. Therefore, $\si$ is homotopic
to an isometric reflection $\si'\in F\acts M'$, whose fixed-point
set is a hypersurface $S'$ homotopic to $S$.  Then $S'$ is a
hyperbolic $n-1$-dimensional manifold homotopy-equivalent to
$\partial N$. Since $\partial N$ is also a hyperbolic manifold, it
follows from Mostow Rigidity Theorem that $\partial N$ and $S'$
are isometric. Let $L'\subset M'$ denote the submanifold bounded
by $S'$ and homotopic to $N'$.

Then we can glue $N$ and $L'$ along their totally-geodesic
boundaries via the isometry $\partial N\to S'$. The result is a
compact hyperbolic manifold $K$ which is obtained by gluing $2m-1$
submanifolds ${O}_j$, each of which is homotopy-equivalent to
${O}$. Since ${O}$ admits a reflection symmetry $\iota$, it
follows (from Mostow Rigidity Theorem)  that $K$ admits an
isometric dihedral group action
$$
D_{2m-1}\acts K,
$$
whose fundamental domain is a submanifold $O''$ with the dihedral
angle $\frac{2\pi}{2m-1}$ along the fixed-point set $V''$ of
$D_{2m-1}$.

Note that the hyperbolic manifolds with boundary $O$, $O'$ and
$O''$ are homotopy-equivalent to each other, where the
homotopy-equivalences restrict to isometries between their
boundaries. The boundary of each manifold is totally-geodesic away
from an $n-2$-dimensional submanifold $V, V', V''$; the dihedral
angles equal $\frac{\pi}{m}, \frac{\pi}{m-1}$ and
$\frac{2\pi}{2m-1}$ respectively.

We now take $m$ copies of $O$, $m-1$ copies of $O'$ and glue them
together (using isometries of the components of $\partial
O\setminus V, \partial O'\setminus V'$) to form an manifold $Q$.
In the manifold $Q$ the submanifolds $V, V'$ are identified with a
codimension 2 submanifold $U$ (isometric to $V\cong V'$).  The
total dihedral angle along $U$ equals
$$
m \frac{\pi}{m} + (m-1) \frac{\pi}{m-1}= 2\pi.
$$
Thus the manifold $Q$ is hyperbolic.

On the other hand, we can glue together
$$
2m-1= m+ (m-1)
$$
copies of the manifold $O''$ to form a hyperbolic manifold $Q'$.
Thus there exists a homotopy-equivalence $h: Q\to Q'$ which
carries copies of $O, O'$ to copies of $O''$. However, if there is
an isometry $h': Q\to Q'$ homotopic to $h$, then $h'$ would have
to carry a copy of $O$ to a copy of $O''$. The latter is
impossible since these orbifolds have different dihedral angles
along $V$ and $V''$. Contradiction. \qed

\begin{defn}
When $n\ge 4$ (and $m$ is sufficiently large), we will refer to
the manifolds $M'$ as {\em Gromov-Thurston examples}.
\end{defn}

\section{Geometric preliminaries}

\subsection{Projective structures}

Let $X$ be a smooth manifold and $G\acts X$ be a real-analytic Lie group action.  

An $(X,G)$-structure on a manifold $M$ is a 
maximal atlas $A=\{(U_i, \phi_i): i\in I\}$ where $U_i$'s are open subsets of $M$ and 
$\phi_i: U_i\to \phi_i(U_i)\subset X$  are charts, so that the transition maps
$$
\phi_j \circ \phi_i^{-1}
$$
are restrictions of elements of $G$. Every $(X,G)$--structure 
on $M$ determines a pair
$$
(dev, \rho)
$$
where $dev: \tilde{M}\to X$ is a local homeomorphism defined on the universal cover of $M$ 
and $\rho: \pi_1(M)\to G$ is a representation so that $dev$ is $\rho$-equivariant.  
The map $dev$ is called the {\em developing map} and $\rho$ is called the {\em holonomy representation} 
of $A$. Conversely, each pair $(dev, \rho)$, where $\rho$ is a homomorphism 
$\pi_1(M)\to G$ and $dev$ is a $\rho$-equivariant local homeomorphism, determine an $(X,G)$--structure on $M$. 

\begin{rem}
Analogous definitions make sense for orbifolds.
\end{rem}

Clearly, every open subset $\Om\subset X$ has a canonical $(X,G)$--structure $can$ induced from $X$. 
If $\Ga\subset G$ acts properly discontinuously and freely on $\Om$ then $can$ projects 
to the quotient manifold $\Om/\Ga$. 

The most relevant examples of $(X,G)$-structures for this paper are:

1. (Real) projective structures, where $X=\R P^n$, $G=PGL(n+1,\R)$ is the group of projective transformations. 

2. Affine structures, where $X=\R^n$, $G= GL(n+1,\R)\ltimes \R^n$ is the group of affine transformations. 

\medskip
Clearly, every affine structure is also projective. Conversely, given any projective structure 
on $M^n$ there is a canonical affine structure on the appropriate line bundle over $M^n$, induced 
by the tautological line bundle over $\R P^n$.

 We refer the reader to \cite{Goldman(rps), Goldman(1990), 
Choi-Goldman} for the foundational material on real-projective structures.

\subsection{Convex sets}

A subset $K\subset \R^{n+1}$ is called a {\em convex homogeneous cone} if it is convex and is invariant 
under multiplication by positive numbers.

A subset $C\subset \R P^n$ is said to be {\em convex} if it either the entire $\R P^n$ or 
is the image of a convex homogeneous cone $\hat{C}\subset \R^{n+1}\setminus \{0\}$ under the projection 
$\R^{n+1}\setminus \{0\}\to \R P^n$. 

An open subset $C\subset \R P^n$ is convex if and only if either $C=\R P^n$ or there exists a linear 
subspace $\R P^{n-1}\subset \R P^n$ such that $C$ is a convex subset of the affine 
space $\A^n=\R P^n\setminus \R P^{n-1}$. 

Suppose that $\Om$ is an open convex subset of $\R P^n$. Then $\Om$ is said to be {\em strictly convex} 
if its frontier contains no nondegenerate segments.

\medskip 
Given a point $x\in \R^n$ and a set $B\subset \R^n$ let
$\Si=Cone_x(B)$ denote the union of all segments $\ol{xb}, b\in
B$. We will refer to $x$ as the {\em tip} and $B$ as the {\em
base} of this cone.

\begin{lem}
\label{simple} If $B$ is convex then $\Si$ is also convex.
\end{lem}
\proof Let $p, q\in \Si$. Then there exist $a, b\in B$ such that
$p\in \ol{ax}, q\in \ol{bx}$. Thus the segment $\ol{pq}$ is
contained in the planar triangle $\Del(a, b, x)$ with the vertices
$a, b, x$. Since $\ol{ab}\subset B$, it follows that
$$
\ol{pq}\subset \Del(a, b, x)\subset \Si. \qed $$

\medskip
Suppose that $A$ is a projective structure on an $n$-manifold $M$. Let 
$\t{A}$ denote the lift of $A$ to the universal cover $\t{M}$ of $M$.

\begin{defn}
The projective structure $A$ is called {\em convex} if $dev: (\t{M},\t{A})\to \R P^n$ 
is either a 2-fold cover or is an isomorphism onto a convex subset in $\R P^n$. 
\end{defn}

In other words, convex projective structures appear as quotients 
$\Om/\Ga$, where $\Om\subset \R P^n$ is convex and $\Ga$ is a properly discontinuous group 
of projective transformations of $\Om$.

\section{A convexity theorem}
\label{convexitythm}

In geometry one frequently constructs geometric objects by gluing together other geometric 
objects. For instance, given hyperbolic $n$-manifolds $M_1, M_2$ with totally geodesic boundary   
and an isometry $\phi: \partial M_1\to \partial M_2$, one constructs a new hyperbolic 
manifold $M=M_1\cup_\phi M_2$ by gluing $M_1$ and $M_2$ via $\phi$. Under some mild assumptions, if 
$M_1, M_2$ are both complete, then so is $M$. (For instance, it suffices to assume 
that the boundaries of both $M_1, M_2$ have positive normal injectivity radius.) Another instance 
of this phenomenon is Poincare's fundamental domain theorem, where instead of gluing manifolds with boundary 
one glues manifolds with corners. 

Recall that in a complete connected Riemannian 
manifold any two points can be connected by a geodesic. Therefore, 
the most natural generalization of the notion of completeness in the category of projective 
structures is {\em convexity}. The problem however is that typically, union of convex sets is 
not convex. Therefore we have to impose further restrictions in order to get convexity. 

Below is a simple example (which I owe to Yves Benoist) 
of failure of convexity of an affine structure built out of convex fundamental domains. 

Let $P$ denote the convex 2-dimensional polygon in $\R^2$ with the vertices 
$$
(1,0), (2,0), (0,1), (0,2). 
$$
Let $A(x)=2x$ and $B$ be the rotation by the angle $\pi/4$. By gluing the sides of $P$ via 
 $A$ and $B$ we obtain an affine structure on the torus. However this structure is not convex since 
 the image of the developing map is $\R^2\setminus \{0\}$.

\medskip 
The main result of this section is a version of Poincare's fundamental domain theorem 
in the context of convex projective structures. We will show that, under some conditions, 
an affine manifold 
obtained by linear gluing of convex homogeneous cones (with infinitely many faces) is again convex. 
Projectivizing this statement we get a similar result for projective structures.

\medskip 
Throughout this section we will assume that $C$ is an open convex
homogeneous cone in $\R^n$ which is different from $\R^n$ itself.

\begin{dfn}
1. An (open) {\em facet} of $C$ is an open convex homogeneous $(n-1)$-dimensional cone
contained in  the boundary of $C$.

2. A {\em codimension $k$} face of $C$ is an open convex
$(n-k)$-dimensional cone contained in  the boundary of $C$. (We
will mostly need this for $k=2$.)
\end{dfn}

We will use the notations $\bar{F}, \bar{C}$, etc. to denote the
closures of faces, cones, etc. Accordingly, we will refer to
closed facets, closed codimension $k$ faces of $C$, etc.

For each face $F$ of $C$ let $Span(F)$ denote the hyperplane in $\R^n$ spanned by $F$.

Let $D$ be a convex subset of $\R^n$, so that $C\subset D\subset
\bar{C}$, and which is obtained by adding to $C$ some of the
 faces.

\begin{lem}
\label{strict}
 Let $x, y\in \bar{C}$ be such that $y\in D$.
 Then the half-open interval $(x, y]:=\ol{x y}\setminus \{x\}$ is contained in $D$.
\end{lem}
\proof The point $y$ belongs to a certain open convex cone $E$
which is either $C$, or  a face of $C$. Let $B\subset E$ be an
open round ball centered at $y$. Let $\Sigma=Cone_x(B)$ denote the
cone with the tip at $x$ and base $B$, which is the union of
segments connecting $x$ to the points of $B$. Then $\Sigma$ has
nonempty interior $\Sigma^0\subset Span(E)$ containing the open
segment $(x, y]$. By convexity, $\Sigma\subset \bar{E}$, hence
$$
(x, y]\subset \Sigma^0\cup \{y\}\subset E\subset C. \qed
$$

Let $\F$ be a certain collection of faces of $F$, so that $C\in \F$. We define a new convex cone
$C'$ as
$$
C'=\bigcup_{F\in \F} F.
$$
Then $C\subset C'\subset
\bar{C}$. For a face $F\in \F$ we let $F'$ denote the closure of $F$ in
$C'$. By abusing the notation we will continue to refer to the
$F'$'s as {\em faces} of $C'$.

\begin{assumption}
\label{ass1} We assume that $C$ and $\F$ are such that:

\begin{enumerate}

\item Each face $F\in \F$ of codimension $\ge 2$ is the
intersection of the higher-dimen\-sional faces $G_i'$.

\item For each pair of distinct facets $F_1, F_2\in \F$, with the
$n-2$-dimensional intersection $\bar{F}_1\cap \bar{F}_2$, we require
this intersection to contain a codimension 2 face $F\in \F$.

\end{enumerate}
\end{assumption}

Suppose now that $X$ is a simply-connected  affine $n$-manifold
obtained by gluing infinitely many copies $C'_j, j\in J,$ of the
cone $C'$ via linear isomorphisms of faces $F'$ of $C'$. We will
refer to the cones $C_j$ as {\em cells} in $X$. We will assume
that each face of $X$ is contained in only finitely many other
faces.

In a similar fashion we define a cell complex $\bar{X}$ by
extending the above gluing maps to the closed cells $\bar{C}_j$.
The {\em origin} $0\in \bar{X}$ is the point corresponding to $0$
in the closed cone $C$. Since $X$ contains infinitely many cells,
the space $\bar{X}$ is not locally compact (the origin is incident
to infinitely many cells).

\begin{rem}
More generally, one can allow spaces $X$ built out of
non-isomorphic convex cones $C_j$. However we do not need this for
the purposes of this paper.
\end{rem}

We then have a {\em developing map} $dev: X\to \R^n$, which is a linear isomorphism on each cell $C_j$. The
developing map extends naturally to a map $dev: \bar{X}\to \R^n$.
This map determines the notion of a {\em segment} $[xy]$ in $\bar{X}$, which is
defined as a path which is mapped by $dev$ homeomorphically to a straight-line segment in $\R^n$.

\begin{defn}
A subset $S$ of $\bar{X}$ is called {\em convex} if every two points in $S$ can be connected by a segment which
is contained in $S$.
\end{defn}

\begin{figure}[tbh]
\centerline{\epsfxsize=2in \epsfbox{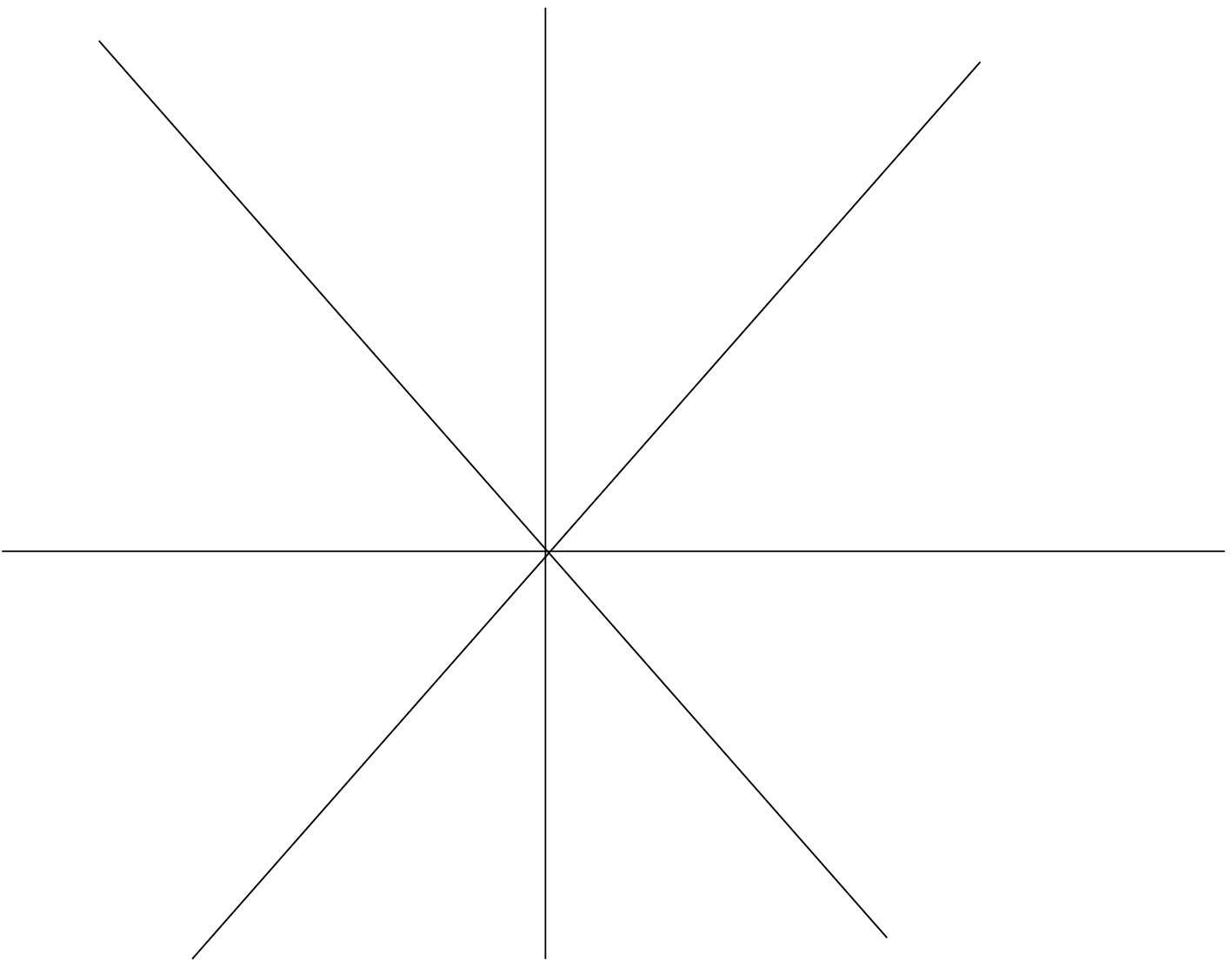}}
\caption{\sl In this example $t=4$.}
\label{f3.fig}
\end{figure}

\begin{assumption}\label{ass2}
1. We assume that $X$ is such that for each point in every
codimension 2 face $E$ of $C_j\subset X$, the tessellation of $X$
by  the adjacent cells is locally isomorphic to a tessellation of
$\R^n$ by cones cut off by a family of $t\ge 4$ hyperplanes
passing through a codimension 2 subspace. (Thus the number of
adjacent cones is $2t\ge 8$.) See Figure \ref{f3.fig}.

2. In addition, we assume that for every pair of cells $C'_i,
C'_j$ sharing a facet $F$, the union $C_i'\cup C_j'$ is convex.
\end{assumption}

Note that this assumption is satisfied in a number of important cases, e.g.
for tessellations corresponding to linear reflection groups.


A {\em wall} $H$ in $X$ is a maximal connected subset of $X$ which
is the union of facets $F_j'$ so that each point $x\in H$ has a
neighborhood $U\subset H$ which is mapped by $dev$
homeomorphically to an open disk (or a half-disk) of a hyperplane
in $\R^n$. Thus the developing map sends each wall to a subset of
a hyperplane in $\R^n$. Therefore, each segment $\si\subset X$ can
intersect a wall $H$ transversally in at most one point.

\begin{rem}
One can show that each wall is a manifold without boundary.
\end{rem}


The main result of this section is the following

\begin{thm}
\label{convex} Suppose that $X$ is as above and $n=dim(X)\ge 3$.
Then:

1. $\bar{X}$ is convex and $dev: \bar{X}\to \R^n$ is an isomorphism onto a convex homogeneous cone in $\R^n$.

2. $X$ is  convex and the developing map $dev$ is an isomorphism
of $X$ onto a proper open convex homogeneous cone in $\R^n$.
\end{thm}

\proof Our proof is modelled on the standard arguments appearing
in the proofs of Poincar\'e's fundamental domain theorem (cf.
\cite{Vinberg(1971)}).

Let $Z$ denote the nerve of the collection of codimension $\ge 2$
faces in $X$. Then $dim(Z)\le 2$; our assumptions on $X$ imply
that $Z$ is a simply-connected regular cell complex. Nonempty
intersections between cells in $Z$ are again cells. Each 2-face
$c$ of $Z$ is a $2t$-gon for a certain $t\ge 4$ (depending on
$c$). Since $C$ and its facets are convex it follows that the
links of vertices of $Z$ contain no bigons. Therefore $Z$
satisfies the {\em small cancellation condition} $C'(1/7)$, see
\cite[Appendix]{Ghys-Harpe(1990)}.

\begin{rem}
Since $Z$ is simply-connected and satisfies $C'(1/7)$ condition,
it follows that $Z$ is contractible. Therefore $\F$ contains only
faces of codimension $\ge 2$.
\end{rem}

However in general $Z$ is not locally compact, since $C$ can have
infinitely many facets. This is where we deviate from Vinberg's
argument \cite{Vinberg(1971)}.

\begin{defn}
A path $p$ in the 1-skeleton $T:=Z^{(1)}$ of $Z$ is a {\em local
geodesic} (or {\em Dehn-reduced}) if it contains no backtracks and no subpaths of length
$>t$ contained in a single $t$-gonal 2-cell of $Z$.
\end{defn}

A {\em degenerate bigon} in $T$ has two vertices $x, y$ and two
{\em equal} edges $\al, \be$, which are local geodesics in $T$.

Given two vertices $x, y\in T$ which belong to a common $t$-gonal 2-cell $c$
of $Z$ and which are distance $t$ apart, there are exactly two
geodesics $\al, \be$ (of length $t$) connecting $x$ to $y$.
The union of these geodesics is the boundary of $c$. We then obtain an {\em elementary bigon}
in $T$ with the vertices $x, y$ and edges $\al, \be$.

More generally, define a {\em corridor} $D$ in $Z$ as a union of 2-cells $F_1,...,F_l$ ($l\ge 2$) of the complex $Z$
so that for each $i$:

1.  $F_i, F_{i+1}$ share an edge $e_i$, called an {\em interior edge} of the corridor $D$.

2. The edges $e_{i-1}, e_{i}$ are ``antipodal'' on the boundary of $F_i$.

Thus each corridor $D$ yields a wall in $X$.

Consider the topological circle $\la$ which is the boundary of the corridor $D=F_1\cup ... \cup F_l$.
Pick two vertices $x, y\in \la$. They are connected by two curves $\be, \ga\subset \la$ whose union is $\la$.

\begin{defn}
In case both $\be, \ga$ are local geodesics, we refer to $\be\cup \ga$ as a {\em simple bigon} with the
vertices $x, y$.
\end{defn}

One can actually show that in case when $\be\cup \ga$ is a simple
bigon, at least one of the arcs $\be, \ga$ is a geodesic in $T$.
Moreover, given a corridor $D$, there are exactly 4 simple bigons
contained in the boundary of the disk $D$.


\begin{lem}
\label{L1} Let $x, y\in T$, be vertices within distance $d$ from
each other. Then:

1. There are only finitely many local geodesics in $T$ connecting
$x$ to $y$.

2. The union $geo(x,y)$ of  geodesics connecting $x$ to $y$ is
convex in $T$. The distance between any two vertices in $geo(x,y)$
is $\le d$. Thus $diam(geo(x,y))=d$.
\end{lem}
\proof Sinc $Z$ satisfies the condition $C'(1/7)$, according to
\cite[Proposition 39, Part (i)]{Ghys-Harpe(1990)}, all (local)
geodesic bigons in $T$ are concatenations of:

1. Degenerate bigons.

2. Elementary bigons.

3. Simple bigons.

\begin{figure}[tbh]
\centerline{\epsfxsize=4.5in \epsfbox{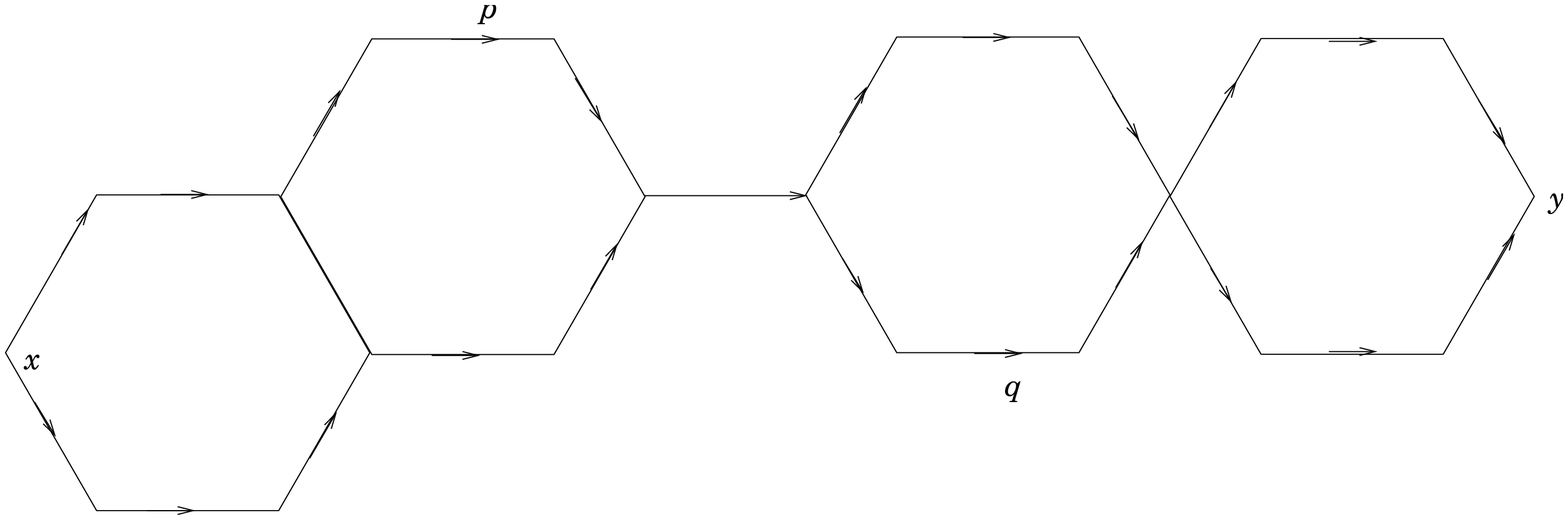}}
\caption{\sl A geodesic bigon with the vertices $x, y$.}
\label{f2.fig}
\end{figure}

See Figure \ref{f2.fig}.

\begin{rem}
The proofs in \cite[Appendix]{Ghys-Harpe(1990)} are given under the assumption that the cell complex
is the Cayley complex of a finitely-generated group.  However the proofs needed for
\cite[Proposition 39, Part (i)]{Ghys-Harpe(1990)} do not require this assumption.
\end{rem}

Consider a local geodesic segment $\ga$ in $T$ with the end-points
$x, y$. Then each subpath $\al$ of length $2$ in $\ga$ is
contained in at most one 2-cell $c_\al$ of $Z$. Thus the union
$$
U:=\bigcup_{\al\subset \ga} c_\al
$$
is compact. According to the above description of geodesic bigons, each geodesic $\be\subset T$
connecting $x$ to $y$ is contained in $U$. Thus the set of such geodesics is finite. This proves the
first assertion.

Each geodesic bigon $\be\cup \ga$ in $T$ bounds a disk $D\subset
Z$ of the least combinatorial area; we will refer to this area as
the {\em area of the bigon}. Then for the given geodesic $\be$
there exists a unique geodesic $\ga$, connecting the end-points
$x, y$ of $\be$, so that the resulting bigon has maximal area.
Namely, consider subpaths $\al\subset \ga$ each of which are
contained in the boundary of a 2-cell $c_\al\subset Z$ and the
length of $\al$ is half of the perimeter of $c_\al$. Then the
distinct subpaths $\al, \al'$ either do not overlap or overlap
along an edge (in the case they belong to a common corridor). Let
$D$ be the union of $\ga$ with all the 2-cells $c_\al$. Then $D$
is the minimal disk bounding the maximal bigon $\be\cup \ga$ with
the vertices $x, y$.

It follows from the description of $D$ that  all  geodesics in $T$
connecting $x$ to $y$ are contained in the 1-skeleton of $D$. If
$z, w\in D^{(0)}$ are not in the same 2-cell contained in $D$,
then each geodesic connecting $p$ to $q$ extends to a geodesic
connecting $x$ to $y$. If $z, w$ belong to the same 2-cell $c$
then all geodesics connecting $z$ to $w$ are contained in
$\partial c$. Thus $D$ is convex. \qed

\medskip
We call a segment $[x,y]\subset \bar{X}$ {\em generic} if the open
segment $(x,y)$ is entirely contained in the union of open cells
and open facets and is not contained in a single facet.

Each generic segment $\si=[x,y]$ determines a path $p(\si)$ in
$T$: The vertices $z_0, z_1, ...,z_k$ and edges of this path correspond to the open
cells $C_0,...,C_k$ and open facets in $X$ which cover the open segment $(x,y)$.

\begin{defn}
\label{gallery}
The union
$$
C_0'\cup...\cup C_k'$$
is the {\em gallery} $gal(p(\si))$ corresponding to the path $p(\si)$. The number $k$ is the
{\em length of the gallery}, it equals the length of the path $p(\si)$
in $T$. By abusing notation we will refer to $k$ as the {\em distance} between $C_0$ and $C_k$.
\end{defn}

The following lemma shows that the distance between $C_0$ and
$C_k$ is the distance between the end-points of the path $p(\si)$
in the graph $T$, thereby justifying the above definition.

\begin{lem}
\label{L2} The path $p(\si)$ is a geodesic in $T$.
\end{lem}
\proof It follows immediately from convexity of the cells in $X$
and the Assumption \ref{ass2}, that $p(\si)$ is a local geodesic.
Suppose that $p(\si)$ is not a global geodesic. Then there exists
a shorter geodesic $q\subset T$ connecting the end-points of
$p(\si)$. Thus $p(\si)\cup q$ is a (locally) geodesic bigon in $T$
bounding a minimal disk $D\subset Z$. Since $p(\si)$ is not a
geodesic, the disk $D$ contains a corridor. Therefore it suffices
to obtain a contradiction in the case when $D$ is a corridor
itself. In this case there is a wall $H$ in $X$ whose intersection
with the gallery $gal(p(\si))$ is not connected. (The wall $H$
passes through the facets in $X$ corresponding to the interior
edges of the corridor $D$.) However the segment $\si$ can
intersect the wall $H$ transversally in at most one point.
Contradiction. \qed

\begin{rem}
The above lemma is analogous to the familiar description of geodesics in Cayley graphs of Coxeter groups.
\end{rem}

We now extend the definition \ref{gallery} to allow galleries
associated with arbitrary geodesic paths $q\subset T$.

Given two cells $A, B$ in $X$, let $Gal(A, B)$ denote the union of
all galleries
$$
A'=C_0'\cup C_1'\cup...\cup C_l'=B',
$$
connecting $A$ to $B$. We define the union of {\em closed
galleries}
$$
Gal(\bar{A}, \bar{B})
$$
as the closure of $Gal(A,B)$ in $\bar{X}$. 
Observe that, according to Lemmata \ref{L1} and \ref{L2}, $Gal(A,
B)$ is a finite union of cells. 
We will prove Theorem \ref{convex} on existence of a segment connecting points $x, y$
by induction on the distance between the cells.

In case when we have points $x, y\in \bar{X}$ (resp. $X$) which belong to the same cell,
there is nothing to prove (since each cell is convex).

Suppose that for every pair of points in $\bar{X}$ which belong to
cells within distance $\le k-1$, there exists a segment
$[x,y]\subset \bar{X}$. Our goal is to prove the same assertion
for $k$.


Pick two cells $A, B$ which are distance $k$ apart and which
correspond to vertices $a, b\in T$. Our goal is to show that there
exists a segment $[x,y]\subset \bar{X}$ for all $x\in \bar{A},
y\in \bar{B}$.

Recall that the union $geo(a,b)$ of geodesics in $T$ connecting
the vertices $a, b$ is convex and has diameter equal to
$d(a,b)=k$. Therefore, each cell $D\subset Gal(A, B)$ which is
adjacent to $B$, the gallery $Gal(A,D)$ is contained in $Gal(A,B)$
and its projection to $T$  has diameter $k-1$. Thus, by the
induction hypothesis, for every pair of points $x', y'\in
Gal(\bar{A},\bar{D})$ the segment $[x',y']\subset \bar{X}$ exists;
moreover, according to Lemma \ref{L2}, this segment is contained
in $Gal(\bar{A},\bar{D})$. Thus $Gal(\bar{A},\bar{D})$ is convex.

Let $\Phi:=\{F_1,..., F_l\}$ denote the set of facets of $B$ which
are contained in the interior of $Gal(A,B)$. Thus for each facet
$F\in \Phi$ there exists a gallery
$$
A'=C_0'\cup C_1'\cup... \cup C_{k-1}' \cup C_k'=B',
$$
so that $F=C_{k-1}'\cap C_k'$.

Fix a point $x$ in the open cell $A$ and let $y\in \bar{B}$ vary.
Let $Y\subset \bar{C}_k$ denote the set of points $y\in \bar{C}_k$
such that there exists a segment $[x,y]\subset \bar{X}$. Each
facet $F\in \Phi$, is contained in a cell $D_F=\bar{C}_{k-1}$
which is distance $k-1$ away from $\bar{C}_0$. Therefore
$$
F\subset Gal(\bar{A}, \bar{D}_F),
$$
and the latter is convex by the induction assumption. Hence $F\subset Y$, which implies that $Y$ is nonempty.

\medskip
Let $B_{sing}\subset B$ denote the (possibly empty) set of points
$y$ such that the segment $[x,y]\subset \bar{X}$ passes through
the origin or through the boundary of a codimension 2 face. Since
$B$ has dimension $\ge 3$ it follows (from the dimension count)
that $L$ does not locally separate $B$. Set $B_{reg}:= B\setminus
B_{sing}$.

\begin{lem}
1. $cl(Y_{gen})\subset Y$.

2. $cl_{B_{reg}}(Y_{gen}) \subset int(Y)\cap B_{reg}$.

3. For each $F\in \Phi$, $F$ is contained in the interior of $Y$.

4. $int(Y)\cap B_{reg}\subset cl_{B_{reg}}(Y_{gen})$.
\end{lem}
\proof 1. Consider a sequence $y_j\in Y_{gen}$ which converges to some $y\in \bar{B}$. Then
$$
[x,y_j]\subset Gal(\bar{A}, \bar{B}).
$$
Since the above union of galleries is compact, it follows that the
sequence of segments $[x,y_j]$ subconverges to a segment $[x,y]$
in $\bar{X}$. Thus $y\in Y$.

2. For each $y\in Y\setminus B_{sing}$ consider the segment
$\si=[x,y]$ and define the point $z=z_{\si}\in [x,y]$ so that
$$
[z,y]\subset \bar{B}, \quad [x,z)\cap \bar{B}=\emptyset.
$$
Consider a sequence $y_j\in Y_{gen}$ which converges to a point
$y\in B$. Then (similarly to the proof of 1), we can assume that
the segments $\si_j=[x,y_j]$ converge to a segment $[x,y]$. For
each $j$ take the point $z_j:= z_{\si_j}$; the limit of these
points is some $z\in \partial B\cap [x,y]$. Without loss of
generality we can assume that all $z_j$'s belong to a common facet
$F$ of $B$, so that $\bar{F}=\bar{D}\cap \bar{B}$, where $D\subset
Gal(A,B)$. Thus $z\in \bar{F}\subset \bar{D}$.

Since $y\in B_{reg}$, there are only three possibilities:

\no (a) $z\in F$.

\no (b) $z$ belongs to an (open) codimension 2 face contained in
$\partial F$.

\no (c) $z$ belongs to the intersection $\bar{D}\cap \bar{B}$ and is not contained in any other closed cell.

In all three cases, by Assumption \ref{ass2}, there exists a
convex neighborhood $U$ of the point $z$ in
$$
Gal(\bar{A}, \bar{D})\cup \bar{B}
$$
so that $U\cap X$ is also convex.
(In case (a) it follows since $X$ is an affine manifold; in case
(b) it follows from Part 1 of Assumption \ref{ass2}; in case (c)
it follows from Part 2 of this assumption.)

Recall that $x\in A \subset int (Gal(\bar{A}, \bar{D}))$, the
latter is convex. Thus, by Lemma \ref{strict},
$$
[x,z)\subset int (Gal(\bar{A}, \bar{D}))\subset X.
$$
By convexity of $\bar{B}$ (and Lemma \ref{strict}), $(z,y]\subset B\subset X$. Pick points
$z'\in int(U)\cap [x,z), z''\in int(U)\cap (z, y]$. Convexity of $X\cap U$ then implies that $z\in X$; hence the entire segment
$[x,y]$ is contained in $X$. Therefore, since $X$ is an affine manifold, there exists a neighborhood $V$ of $y$ in $B$
such that for each $y'\in V$, there exists a segment $[x,y']\subset X$. Thus $y\in int(Y)$.

3. The proof of this assertion is analogous to the last part of the proof of 2: For each $y\in F$ the segment $[x,y]$ is contained
in $X$ and thus we can use the neighborhood $V$ of $y$ as above.

4. For $y\in int(Y)$, let $V\subset Y$ be an open ball containing
$y$. Since the union of facets of the cell $C$ is dense in
$C'\setminus C$, we see that $V\setminus Y_{gen}$ is nowhere dense
in $V$. Hence $Y_{gen}\cap V$ is dense in $V$ and therefore $y\in
cl(Y_{gen})$. The assertion 4 follows. \qed

\begin{cor}
$\bar{B}=Y=cl(Y_{gen})$.
\end{cor}
\proof First, by combining 2 and 4 we see that
$$
cl_{B_{reg}}(int Y \cap B_{reg})= cl_{B_{reg}}(Y_{gen}).
$$
This of course implies that $cl(int Y \setminus
B_{sing})=cl(Y_{gen})$.

Thus $int (Y\cap B_{reg})$ is both closed and open in $B_{reg}$.
By 3, the set $int (Y\cap B_{reg})$ is nonempty. Since $B_{reg}$
is connected (recall that $B_{sing}$ does not locally separate),
we conclude that $int (Y\cap B_{reg})=B_{reg}$. Thus
$B_{reg}\subset Y$. Since $B_{reg}$ is dense is $\bar{B}$, by
applying 1 we see that
$$
\bar{B}\subset cl(B_{reg}) = cl( int (Y\cap B_{reg}) )=
cl(Y_{gen})\subset Y.
$$
Thus $\bar{B}=Y=cl(Y_{gen})$. \qed

Therefore, for each point $y\in \bar{B}$ there exists a segment $[x,y]\subset \bar{X}$. Recall that we assumed that $x\in A$.
For a point $y\in \bar{B}$ and $x\in \bar{A}$, pick a sequence $x_j\in A$. Then there exist a sequence
$y_j\in Y_{gen}$ which converges to $y$; the segments $[x_j, y_j]$ are all contained in $Gal(A, B)$. Therefore, by
compactness of  $Gal(\bar{A}, \bar{B})$, we conclude that the segments $[x_j, y_j]$ converge to a segment $[x,y]\subset \bar{X}$.

Thus we proved that for each pair of cells $A$ and $B$ within distance $k$, and each pair of points $x\in \bar{A}, y\in \bar{B}$,
there exists a segment $[x,y]\subset \bar{X}$.

\begin{rem}
Moreover, $[x,y]\subset Gal(\bar{A}, \bar{B})$.
\end{rem}

Hence, by induction, we conclude that $\bar{X}$ is convex. It follows that $dev: \bar{X}\to \R^n$ is a continuous bijection
onto a convex cone $K\subset \R^n$. This proves the first assertion of Theorem \ref{convex}.


Our next goal is to prove convexity of $X$. Let $x, y\in X$. Pick a relatively compact neighborhood $U$ of $x$ in $X$.
Then $U$ is covered by finitely many cells $C_j'\subset X$. Choose a cell $A'$ containing $y$.
For each $x'\in U$ the segment $[x',y]\subset \bar{X}$ is covered by a finite union of faces contained in
$$
\bigcup_{j} Gal(\bar{C}_j, \bar{A}).
$$
Therefore the cone $Cone_y(U)$ with the tip $y$ and the base $U$
is covered by finitely many closed cells. Thus the developing map
$dev$ sends $Cone_y(U)$ homeomorphically onto a convex subset
$$
\Si:=dev(Cone_y(U))=Cone_{dev(y)}( dev(U))\subset \R^n.
$$
Clearly, the open segment $(dev(x), dev(y))$ is contained in the interior of $\Si$.
It follows that the open segment $(x,y)$ is also contained in the interior of
$$
\bigcup_{j} Gal(\bar{C}_j, \bar{A}).
$$
Hence the open segment $(x,y)$ is contained in $X$. Thus $dev$ is
a homeomorphism of $X$ onto a convex homogeneous cone in $\R^n$.
Properness of this cone follows from infiniteness of the number of
cells in $X$. \qed

\section{Products of matrices}
\label{pro}

 In this section we will consider the following
problem:

\begin{prob}
Let $G$ be a Lie group with a collection of 1-parameter subgroups
$G_1,...,G_k\subset G$. Analyze the image of the map
$$
Prod: \prod_{i=1}^k G_i \to G
$$
given by $Prod(g_1,...,g_k)=g_1\cdot...\cdot g_k$.
\end{prob}

In the case when $G=SO(3)$, this problem is ultimately related to
the  variety of geodesic $k$-gons  in $S^3$ with the fixed
side-lengths, \cite{Kapovich-Millson1997}.
(See \cite{Kapovich-Millson(1992)}, \cite{Maubon},
for the relation of this product problem to bending deformations of
flat conformal structures.)

Here we consider the case of $G=GL(2,\R)$; the
subgroups $G_i$ are orthogonal conjugates of the group of diagonal
matrices $\{Diag(1, e^t), t\in \R\}$. More specific problem then
is:

\begin{prob}
\label{P2}
Show that under appropriate conditions on the subgroups $G_i$, the image of the map
$Prod$ contains the subgroup $SO(2)\subset GL(2,\R)$.
\end{prob}

Let ${\mathfrak g}{\mathfrak l}(2,\R)= {\mathfrak p}\oplus
{\mathfrak o}(2)$ denote the Cartan decomposition of the Lie
algebra of $GL(2,\R)$. The Lie algebras ${\mathfrak p}_i$ of
$G_i$'s are contained in ${\mathfrak p}$. Let $e:=(1,...,1)\in
\prod_{i} G_i$. Then derivative
$$
d Prod_{e}: \oplus_{i} {\mathfrak p}_i \to {\mathfrak g}{\mathfrak
l}(2,\R)
$$
is the map
$$
(\xi_1,...,\xi_k)\mapsto \sum_{i} \xi_i.
$$
Therefore its image is contained in ${\mathfrak p}$ and hence is
orthogonal to ${\mathfrak o}(2)$. Thus one cannot approach Problem
\ref{P2} by making infinitesimal calculations.

There is probably a purely algebraic or analytic solution to Problem
\ref{P2}; we will use hyperbolic geometry instead. Given a basis
$(v, w)$ of $\R^2$ and $t\in \R$ we define the matrix
$$
A=A_{v,w,t}
$$
to be the linear transformation which fixes $v$ and sends $w$ to $e^t w$.

Consider the projective action of $GL(2,\R)$ on the circle $\R
P^1$ (which we identify with the boundary of the hyperbolic plane
$\H^2$). We will use the notation $[A]\in PGL(2,\R)$ for the
projection of the matrix $A\in GL(2,\R)$.

The vectors $v, w$ project to fixed points $[v], [w]$ of the
projective transformation $[A]$. We identify the hyperbolic plane
$\H^2$ with the unit disk in $\R^2$ in such a way that the group
$O(2)\subset GL(2,\R)$ fixes the origin $0$ in $\H^2$. Then the
hyperbolic geodesic $L_A=\ol{[v] [w]}\subset \H^2$ invariant under
$[A]$ passes through the origin $0$ (and hence is a Euclidean
straight line). We parameterize the geodesic $L$ with the unit
speed and orient $L$ in the direction from $[v]$ to $[w]$, thereby
identifying it with the real line. The origin in $\H^2$
corresponds to zero in $\R$. We let $L^{\pm}$ denote the positive
and negative rays (starting at $0$) in $L$ corresponding to this
orientation.

In these coordinates, the isometry $[A]$ acts on $L$ by $r\mapsto
r+t$. The isometry $[A_{w,v,t}]$ acts on $\H^2$ by the translation
$r\mapsto r-t$ along the geodesic $L$. By considering the action
of $SO(2)$ by conjugation we see the following:

Let $R=R_\phi\in SO(2)$ be the rotation by the angle $\phi$. Then the matrix
$$
R_\phi A_{v,w,t} R_{\phi}^{-1}= A_{R(v), R(w),t}
$$
 acts on $\H^2$ by translation $r\mapsto r+t$ along the geodesic
$$
\ol{ R_{\phi/2}([v]) R_{\phi/2}([w])}.
$$

We now assume that we are given 1-parameter groups
$$
G_1= \{ A_{e_1, e_2, t}: t\in \R\}, \quad G_2= R_{\pi/4} G_1 R_{-\pi/4},$$
$$
G_3=R_{\pi/2} G_1 R_{-\pi/2}= \{ A_{e_2, e_1, t}: t\in \R\}, \quad
G_4= R_{-\pi/4} G_1 R_{\pi/4}.
$$
Geometrically, these are groups of translations
along two orthogonal hyperbolic geode\-sics $L_1$ and $L_2$ in $\H^2$ ($G_i$ and $G_{i+2}$ translate along $L_i$
in the opposite directions, $i=1, 2$). Given a matrix $A_i\in G_i$ we let $\ell_i:=\ell(A_i)$ denote the translation
length of $[A_i]$ along its invariant geodesic; here we are ignoring the orientation so that $\ell_i\ge 0$.

Thus, in order for $A_i\in G_i, i=1,...,4$ to have the product equal to
$R_{\phi}$ it is necessary and sufficient to have:

1. The product of the eigenvalues of $A_i$'s is equal to $1$ (i.e.
the product of four matrices is in $SL(2,\R)$). Equivalently,
$$
t_1 +t_2 +t_3 +t_4=0.
$$

2. The product of the hyperbolic translations
$$
[A_4]\circ [A_3]\circ [A_2]\circ [A_1]
$$
is the rotation $R_{\phi/2}$ around the origin in $\H^2$. In particular, the above product of hyperbolic isometries
has to fix the intersection $L_1\cap L_2$.

\begin{rem}
Similar description, of course, will be valid for more general choices of 1-parameter
groups $G_i$ which are conjugate to $G_1$
by rotations $R_{\theta_i}$, $i=1,...,k$.
\end{rem}

With this geometric interpretation it is clear, for instance, that
the product $A_3 A_2 A_1$ is never a nontrivial rotation. The
reason is that unless  $A_2=1$, $[A_1]=[A_3]^{-1}$, the product of
the hyperbolic isometries does not fix the origin.

\medskip
We now make the situation a bit more symmetric and require that
$$
\ell_1=\ell_4, \quad \ell_2=\ell_3.
$$
In particular, the Condition 1 will be satisfied provided that
$$
t_1, t_3>0,\quad t_2, t_4<0.
$$

We then consider the images of the origin under the compositions
of the isometries $[A_1], [A_2], [A_3], [A_4]$. We let $x_0:=0$;
$x_i:= [A_i]( x_{i-1})$, $i=1,...,4$. Given a number $\theta\in
(-\frac{\pi}{2}, 0)$ set
$$
\al=\al(\theta):= \frac{\pi}{2} +\theta.
$$

\begin{lem}
For every $\theta\in (-\frac{\pi}{2}, 0]$ there exists
a pair of continuous functions $\ell_i=\ell_i(\theta), i=1,2,$ so
that:

1.
$$
\left\{ \begin{array}{c} \cosh(\ell_2)=\cosh(\ell_1) \sin(\al)\\
\sinh^2(\ell_1)=\cos(\al) \end{array}\right.
$$
In particular, $\ell_i(0)=0$, $i=1, 2$.

2. For $\ell_i=\ell_i(\theta)$ the composition
$$
[A_4]\circ [A_3]\circ [A_2]\circ [A_1]
$$
is the (counter-clockwise) rotation $R_\theta$ around the origin $0\in \H^2$ by the
angle $\theta$.
\end{lem}
\proof Let $L_1, L_2$ be the pair of oriented geodesics in $\H^2$
(invariant under the subgroups $[G_1]=[G_3]$, $[G_2]=[G_4]$
respectively) which intersect orthogonally at the origin.

We orient the geodesics $L_1, L_2$ away from the points $[e_1]$,
$[R_{\pi/4}(e_1)]$ fixed by $[A_1]\in [G_1], [A_2]\in [G_2]$. Let
$L_i^+$ denote the positive half-rays in these geodesics.

Recall that $\al\in (0, \frac{\pi}{2}]$ and that we will be using the translation parameters so that
$$
t_1, t_3>0,\quad t_2, t_4<0.
$$

The key observation is that there exists a unique geodesic
quadrilateral (a {\em Lambert quadrilateral}) $Q_{\al}=[0, y_1, x_2, y_2]$ in $\H^2$ with the
three right angles (at the vertices $0, y_1\in L_1^+, y_2\in
L_2^+$) and the angle $\al$ at the vertex $x_2$.  See Figure
\ref{f1.fig}. The orientation on $Q_{\al}$ given by the ordering
of its vertices is clockwise, which corresponds to the assumption
that $\theta\le 0$.

Set
$$
\ell_2=\ell_3:=d(0, y_1)=d(0 y_2)
$$
and
$$
\ell_1=\ell_4:= d(y_2, x_2)= d(y_1, x_2).
$$
It is clear that $\ell_1, \ell_2$ are continuous functions of
$\theta$ so that $\ell_i(0)=0$. The equations relating $\al,
\ell_1, \ell_2$ follow immediately from the hyperbolic
trigonometry, see \cite[Theorem 7.17.1]{Beardon}.

Choose points $x_1\in L_1^+, x_3\in L_2^+$ so that
$$
d(0, x_1)=\ell_1=d(0, x_3)=\ell_4.
$$
Now take the hyperbolic translations $g_1, g_3$ along $L_1$
sending $0$ to $x_1$ and $y_1$ to $0$ respectively. Define the
hyperbolic translations $g_2, g_4$ along $L_2$ sending $0$ to
$y_2$ and $x_3$ to $0$ respectively. Thus the isometries $g_1,
g_4$ have the translation lengths $\ell_1=\ell_4$; the isometries
$g_2, g_3$ have the translation lengths $\ell_2=\ell_3$.

It is clear from the Figure \ref{f1.fig} that
$$
g_2(x_1)=x_2, g_3(x_2)=x_3
$$
and therefore
$$
g_4\circ g_3\circ g_2\circ g_1(0)=0.
$$

Hence the above composition of translations is a certain rotation
$R_\phi$ around the origin. In order to compute the angle $\phi$
of rotation take two vectors $\xi_1, \xi_2\in T_0\H^2$ tangent to
the geodesic rays $L_i^-$, $i=1,2$. Then the images of $\xi_1,
\xi_2$ under
$$
d(g_2\circ g_1),\quad d(g_3^{-1} \circ g_4^{-1})
$$
are tangent to the geodesic segments $\ol{x_2 y_2},\ol{x_2 y_1}$
respectively. Therefore the angle $\phi$ equals
$\al-\frac{\pi}{2}$ (the rotation is in the clockwise direction).
Thus $\phi=\theta$. \qed

\begin{figure}[tbh]
\begin{center}
\input{f1b.pstex_t}
\end{center}
\caption{\sl  Lambert quadrilateral.}
\label{f1.fig}
\end{figure}

Let $A_i\in G_i$ denote the matrices corresponding to the
hyperbolic translations $g_i$. 

\begin{cor}
$$
A_4\cdot A_3\cdot A_2\cdot A_1=R_{\frac{\al}{2}-\frac{\pi}{4}}.
$$
\end{cor}

Therefore we get the following:

\begin{thm}
\label{product}
For each $\tau\in (-\frac{\pi}{4},0]$ there is a unique pair of numbers $t_1\ge 0, t_2\le 0$ so that for the set of parameters
$$
 \overrightarrow{t}=(t_1, t_2, -t_2, -t_1)
$$
the product of the corresponding matrices equals the rotation $R_{\tau}$.
Moreover, the function $\overrightarrow{t}$ depends continuously
on $\tau$.
\end{thm}

{\bf Projective generalization.}

\medskip
Let $P\subset \R P^n$ be a projective hyperplane, $p\in \R P^n
\setminus P$ and $t\in \R$. Then there exists a unique map
$A=A_{P,p,t}\in PGL(n+1,\R)$ satisfying:

1. $A$ fixes $P\cup \{p\}$ pointwise.

2. The derivative $dA_{p}$ equals $e^t I$.

\medskip
Suppose now that $n=2$, $\R^2$ is the affine patch of $\R P^2$.
Let  $P\subset \R P^2$ be the projective line tangent to the unit
vector $v\in T_0 \R^2$, $p\in \R P^2\setminus \R^2$ be the point
at infinity so that the corresponding line $l$ through the origin
contains the unit vector $w$ orthogonal to $v$. Then
$$
A_{P,p,t}=A_{v,w,t}\in GL(2,\R)\subset PGL(3,\R).
$$
 The identity
extension of this linear transformation to the element $\hat{A}\in
GL(n,\R)\subset PGL(n+1,\R)$ equals
$$
A_{Q,q,t}
$$
where $Q$ is the projective hyperplane through the origin
orthogonal to $w$, the point $q\in \R P^{n}\setminus \R^n$
corresponds to the line $l$ as above.

Consider now a collection $P_1, P_2, P_3, P_4$ of projective
hyperplanes in $\R P^n$ passing through the origin, so that the
intersection
$$
\bigcap_i P_i =S
$$
is a codimension 2 projective hyperplane in $\R P^n$. We assume
that the consecutive hyperplanes intersect at the angles
$\frac{\pi}{4}$. For each $P_i$ let $p_i\in  \R P^n\setminus \R^n$
be the ``dual point'' i.e. the corresponding line $l_i$ through
the origin is orthogonal to $P_i$.

\begin{rem}
Somewhat more invariantly, one can describe this setting as
follows. We fix a positive definite bilinear form on $\R P^n$ so
that the points $p_i$ are dual to the hyperplanes $P_i$. Therefore
the assumption that
$$
\bigcap_i P_i =S
$$
is a codimension 2 projective hyperplane in $\R P^n$ implies that
$\{p_1,...,p_4\}$ is contained in a projective line $s\subset \R
P^n$ dual to $S$. We then are assuming that the points
$p_1,...,p_4$ are cyclically ordered on $s$ so that the distance
between the consecutive points is $\pi/4$. (Note that  $\R P^1=s$
has length $\pi$.)
\end{rem}

Then Theorem \ref{product} can be restated as follows:

\begin{thm}
\label{composition}
 For each angle $\tau\in (-\frac{\pi}{4}, 0]$ there is a unique set of parameters
$\overrightarrow{t}=(t_1, t_2, t_3, t_4)=(t_1, t_2, -t_2, -t_1)$ with $t_1\ge 0, t_2\le 0$
so that the composition of the corresponding projective
transformations
$$
A_{P_4,p_4,t_4}\circ ... A_{P_1,p_1,t_1}
$$
equals the rotation $R_{\tau}$ around $S$ by the angle $\tau$,
fixing $S$ pointwise. Moreover, the function $\overrightarrow{t}$
depends continuously on $\tau$.
\end{thm}

\section{Bending}
\label{bending}

In this section we review the {\em bending deformation} of
projective structures.

Recall that in the end of the previous section we defined
projective transformations
$$
A=A_{P,p,t}\in PGL(n+1,\R)
$$
corresponding to the triples $(P,p,t)$, where $P\subset \R P^n$ is
a projective hyperplane, $p\in \R P^n \setminus P$ and $t\in \R$.

\begin{figure}[tbh]
\begin{center}
\input{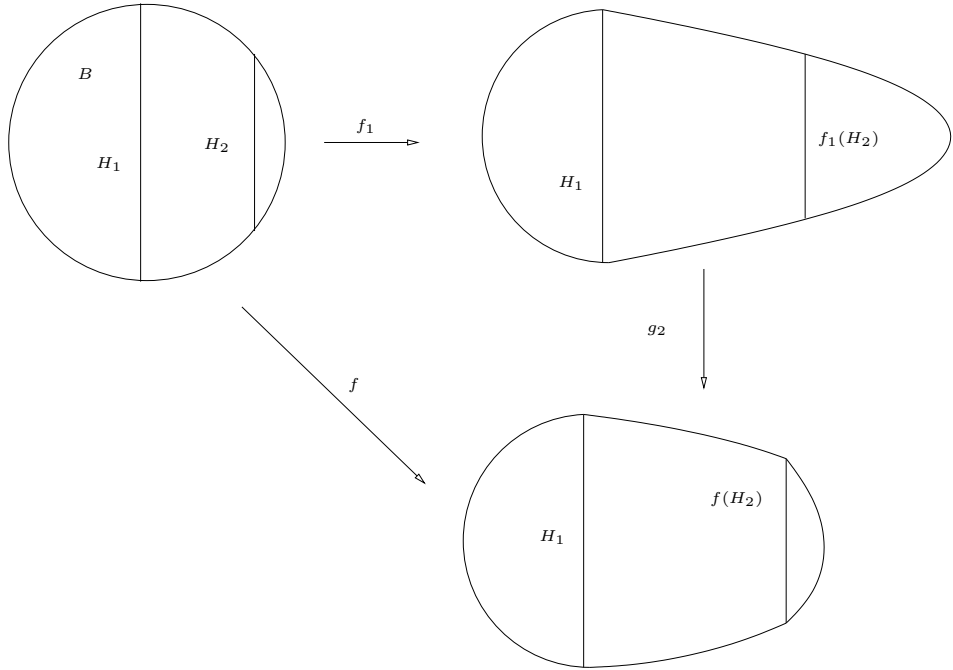}
\end{center}
\caption{\sl Projective bending.}
\label{fig7}
\end{figure}


Before proceeding with the general definition we start with a
basic example of bending. Let $B$ denote the open unit ball in
$\R^n$, which we will identify with the hyperbolic $n$-space. Let
$H_1,...H_k$ denote disjoint hyperbolic hypersurfaces in $B$ so
that $H_i$ separates $H_{i-1}$ from $H_{i+1}$, $i=2,...,k-1$. We
assume that $H_i$'s are cooriented in such a way that $H_{i+1}$ is
to the right from $H_i$, $i=1,...,k$. Let $H_i^{\pm}$ denote the
half-space in $B$ bounded by $H_i$ and lying to the left (resp.
right) from $H_i$. We set
$$
B_i:= H_i^+\cap H_{i+1}^-, \quad i=1,...,k-1$$
and
$$
B_0:= H_1^-, B_k:= H_k^+.
$$

Let $P_i$ denote the projective hyperplane containing $H_i$ and
let $p_i\in \R P^n$ denote the point dual to $P_i$ with respect to
the quadratic form where $B$ is the unit ball. Choose real numbers
$t_1,...,t_k$. Our goal is to {\em bend} $B$ projectively in $\R
P^n$ along the hypersurfaces $H_i$ with the bending parameters
$t_i$, $i=1,...,k$. Let $A_i:= A_{P_i,p_i, t_i}$, $i=1,...,k$.

We will do bending inductively. First, let $f_1: B\to \R P^n$
denote the map which is the identity on $H_1^-$ and $A_1$ on
$H_1^+$. We then would like to bend $B_1:= f_1(B)$ along
$f_1(H_2)$. The corresponding bending map $g_2$ is the identity on
$f_1(H_2^-)$ and $A_2'$ on $f_1(H_2^+)$, where
$$
A_2'= A_1 \circ A_2 \circ A_1^{-1}.
$$
Therefore the map
$$
f_2: B\to g_2(B_1)
$$
equals to $id$ on $B_0=H_1^-$, to $A_1$ on $B_1=H_1^+\cap H_2^-$
and to $A_1\circ A_2$ on $H_2^+$. Continuing in the fashion
inductively we eventually obtain the bending map
$$
f: B\to \R P^n
$$
so that the restriction $f|B_i$ equals
$$
A_1\circ ... \circ A_i, i=1,...,k-1,
$$
and $f|B_0=id$. The same construction works for an arbitrary
locally finite collection ${\mathcal H}$ of disjoint hyperplanes
$H_i$. We then pick a component $C_0$ of $Y=B\setminus \cup_i H_i$
where the bending map $f$ is the identity. Given a component $C_k$
of $Y$ we take the finite subcollection $\{H_1,...,H_k\}$ of
hyperplanes in ${\mathcal H}$ separating $C_0$ from $C_k$. We then
repeat the above construction of bending map to define the
restriction of bending to $C_k$.

\begin{rem}
More generally, if $(L,\mu)$ is a measured codimension 1
totally-geodesic lamination in $B$, we can define projective
bending with respect to this lamination. However, in view of
Ratner's theorem,  this generalization is not useful in the
context of bendings of compact manifolds of dimension $\ge 3$.
\end{rem}

We now give the general definition of bending.

Let $M$ be a projective manifold (or, more generally, an
orbifold). Let $f: \t{M}\to \R P^n$ and $\rho: \Ga=\pi_1(M)\to
PGL(n+1,\R)$ be the developing map and the holonomy of $M$.

Let $L\subset M$  be a proper hypersurface (possibly contained in
the boundary of $M$); let $\t{L}\subset \t{M}$ be the preimage of
$L$ in the universal cover of $M$.

We call the hypersurface $L$ {\em flat} if it satisfies the
following:

\begin{enumerate}
\item  Each point $x\in L$ has a neighborhood $U\subset L$ so that
the developing map sends $U$ to an open subset of a projective
hyperplane in $\R P^n$.

\item For a component $\t{L}_i\subset \t{L}$ let $\Ga_i$ be the
stabilizer of $\t{L}_i$ in $\Ga$. Then the group $\rho(\Ga_i)$
stabilizes a projective hyperplane $P_i=Span(f(\t{L}_i))\subset \R
P^n$.

We then require that for each $\t{L}_i$, the group $\rho(\Ga_i)$
has an isolated fixed point $p_i\in \R P^n$ which is disjoint from
$P_i$.

\end{enumerate}

We define a {\em cooriented  lamination} $L$ in $M$ as follows.
Consider the union $L$ of flat connected cooriented hypersurfaces
$L_i$ in $M$, which intersect the boundary of $M$ transversally
and so that for distinct $L_i, L_j$
$$
L_i\cap L_j \cap int(M)=\emptyset.
$$

In this paper we will be assuming that the collection of
hypersurfaces $L_i$ is locally finite in $M$, although one can
make the discussion more general.

\begin{defn}
A {\em transverse measure} for $L$ is a locally constant function
$\mu: L\cap int(M)\to (0,\infty)$, $\mu: L_i\mapsto e^{t_i}$.
\end{defn}

A  measured cooriented lamination is the pair $\la=(L,\mu)$. The
measured lamination $\la=(L,\mu)$ lifts to a cooriented measured
lamination $\t\la=(\t{L}, \t\mu)$ in $\t{M}$.

We now define the {\em bending deformation} $c_\la$ of the
projective structure $c$ on $M$ along the lamination $\la$. The
structure $c_\la$ will have the developing map $f_\la$ satisfying
the following properties:

For each component $H=\t{L}_i$ of $\t{L}$ with the stabilizer
$\Ga_i$; let $P_i ,p_i$ denote the projective hyperplane and a
point in $\R P^n$ stabilized by $\rho(\Ga_i)$ as above. Let $H_-,
H_+$ denote the components of $\t{M}\setminus \t{L}$ to the left
and to the right of $H$ (with respect to the coorientation). We
then require that there exists a projective transformation $g\in
PGL(n+1,\R)$ so that
$$
f_\la| H_-= g\circ f|H_-,
$$
$$
f_\la|H_+ = A_{P,p,t}\circ g \circ f|H_+.
$$

It is clear that the map $f_\la$ with these properties exists and
is unique up to postcomposition with projective transformations of
$\R P^n$. By construction, $f_\la$ is a local homeomorphism. Since
$\t\la$ is $\Ga$-invariant, it follows that the map $f_\la$ is
equivariant with respect to a homomorphism $\rho_\la: \Ga\to
PGL(n+1,\R)$.

Thus the pair $(f_\la, \rho_\la)$ determines a projective
structure $c_\la$ on $M$.

\bigskip
The following simple lemma is used to ensure convexity of the
projective structures on Gromov-Thurston examples.

Let $H$ be a hyperplane in $\R^n$; let $H_{\pm}$ denote the closed
half-spaces in $\R^n$ bounded by $H$. Suppose that $D\subset H_-$
is a compact convex subset; let $F$ denote the intersection $H\cap
S$. Pick a point $p\in H_+$. Let $\Si=Cone_p(F)$ denote the
(convex) cone with the vertex $p$ and the base $F$.

\begin{lem}\label{un}
Suppose that for each $x\in D$ the segment $\ol{xp}$ crosses $H$
inside $F$. Then the union $D\cup \Si$ is convex.
\end{lem}
\proof Clearly, $D\cup \Si=Cone(p, D)$. Now convexity follows from
Lemma \ref{simple}. \qed

Suppose that $p_i\in H_+$ is a sequence of points so that $D, H,
p_i$ satisfy all the above conditions. Assume that $\lim_i
p_i=p\in \R P^n \setminus \R^n$; let $\Si\subset \R P^n$ denote
the limit of the cones $Cone_{p_i}(F)$. Let $E\subset H_+$ be 
a compact subset contained in the cone $\Si$ so that
$$
D\cap E=D\cap H=E\cap H=F.$$

\begin{cor}
\label{C1}
 $D\cup E$ is convex.
\end{cor}
\proof By taking the limit, Lemma \ref{un} implies that for each
pair of points $x_1\in \Si\setminus \{p\}$, $x_2\in D$, the
intersection $\ol{x_1 x_2}\cap H$ is contained in $F$. Since $E$
is a convex subset of $\Si$; it follows that $\ol{x_1 x_2}\subset
D\cup E$. \qed

Let $P$ denote the projective closure of the hyperplane $H$.

\begin{cor}
\label{C2} For each $t\in \R$, the union
$$
A_{P,p,t}(E)\cup E$$
 is convex.
\end{cor}
\proof The set $E_t=A_{P,p,t}(E)$ is clearly convex. By the definition
of $A_{P,p,t}$, the set $E_t$ is
contained in the cone $\Si$ and $E_t\cap H=E_t\cap D=F$. \qed

\section{Construction of convex projective\\ structures on Gromov-Thurston manifolds}
\label{last}

Assume as before that the hyperbolic manifold $M$ satisfies
Assumption \ref{ass0} and let $M'$ be the manifold constructed
(using pieces of $M$) as in section \ref{gte}.

\begin{thm}
For each natural number $m\ge 8$ which is divisible by $4$, the
manifold $M'$ admits a convex projective structure.
\end{thm}

\proof The proof breaks in two steps:

1. We first {\em bend} the (hyperbolic) projective structure $c$
on the manifold with boundary $N'\subset M'$ (see section
\ref{gte}) in order to obtain a new projective structure $c_\la$
which has {\em flat boundary}. We then construct a projective
structure $a'$ on $M'$ by gluing two copies of $(N', c_\la)$ {\em
via an order 2 rotation}.

2. We use the Theorem \ref{convex} to verify that $(M', a')$ is
convex.

\begin{figure}[tbh]
\begin{center}
\input{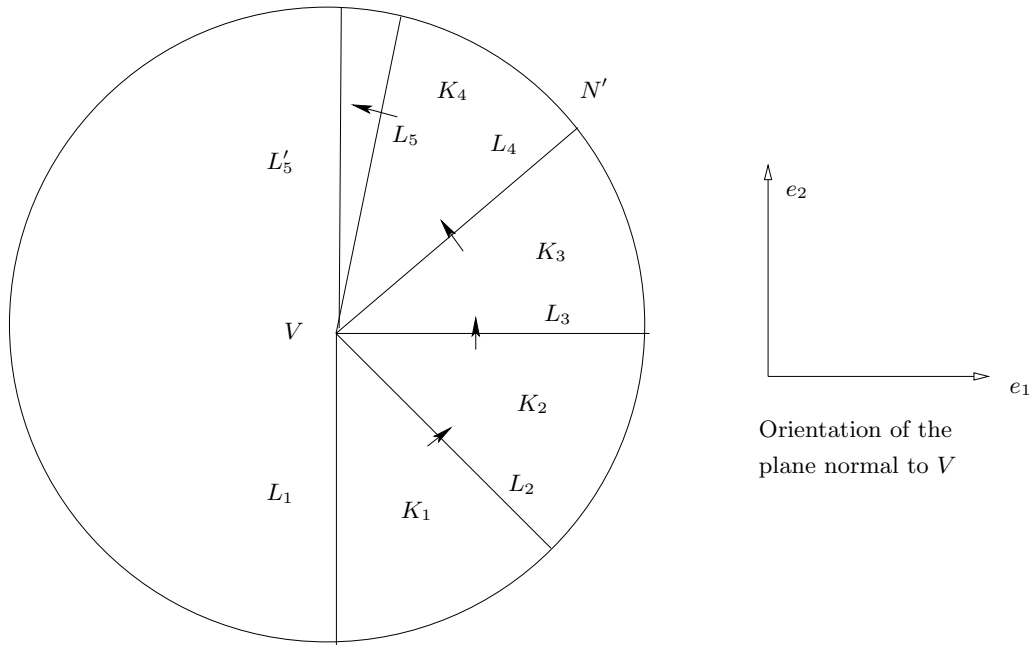}
\end{center}
\caption{\sl Coorientation on $L$. }
\label{fig5}
\end{figure}

\medskip
{\bf Step 1.} We begin by observing that since $m$ is divisible by
$4$, the group $D_m$ contains the dihedral subgroup $D_4$. The
fixed-point sets of the reflections contained in $D_4$ yield
codimension 1 totally-geodesic submanifolds $L_0, ..., L_4\subset
N'$, see figure \ref{fig5}. The angle between $L_i, L_{i+1}$
equals $\frac{\pi}{4}$, $i=0,...,3$.

We assume that the submanifold $N'\subset M$ is chosen in such a
way that $L_0\subset \partial N'$. Then the boundary of $N'$ is
the union $L_0\cup L_0'$.

We let $c$ denote the projective structure on $N'$ defined by the
hyperbolic structure on $N'$. Since each $L_j$ is totally-geodesic
(in the hyperbolic manifold $M$) it is {\em flat} as a
hypersurface in the projective manifold $(N',c)$. We define the
cooriented lamination $L=L_1\cup...\cup L_4$ by coorienting
$L_j$'s as in figure \ref{fig5}. (Each arrow indicates the
direction from left to right.)

We identify the hyperbolic space $\H^n$ with the unit ball $B$ in
the Euclidean space $\R^n\subset \R P^n$, so that each flat
hypersurface $L_0,...,L_4, L_0'$ corresponds (under the developing
map) to a projective hypersurface $P_0,..., P_4, P_0'\subset \R
P^n$ passing through the codimension 2 hyperplane $Q\subset \R
P^n$ containing the origin. The angles between $P_i, P_{i+1}$
equal $\frac{\pi}{4}$ and the angle between $P_0$ and $P_0'$ is
$\pi-\frac{\pi}{m}$.

For each $i=1,...,4$ let $p_i\in \R P^n\setminus \R^n$ be the
point dual to $P_i$ (see section \ref{pro}).

We orient the 2-plane normal to $Q$ in $\R^n$ so that the
orientation agrees with the above coorientation of $L$.

Then, since $m>4$, according to Theorem \ref{pro}, there exist
real numbers $t_1, t_2, t_3=-t_2, t_4=-t_1$ so that
\begin{equation}
\label{1} A_{P_4,p_4,t_4}\circ ... A_{P_1,p_1,t_1}=
R_{-\frac{\pi}{m}},
\end{equation}
is the rotation around $Q$ by the angle $-\frac{\pi}{m}$.

The rotation $R_{-\frac{\pi}{m}}$ sends $P_0'$ to $P_4$. Let
$\mu$ denote the transverse measure to $L$ defined by $L_i\mapsto
e^{t_i}$. Let $\la:=(L,\mu)$ and $c_\la$ be the projective
structure on $N'$ obtained from $c$ by bending along $\la$.

The equation (\ref{1}) implies that the projective manifold $(N',
c_\la)$ has {\em flat boundary}. The boundary manifold is actually
hyperbolic and is isometric to the boundary of the hyperbolic
manifold $N$ (which has geodesic boundary). Let $\theta':
(\partial N', c_\la) \to (\partial N', c_\la)$ denote the
isometric involution which interchanges $L_0$ and $L_0'$ (see
section \ref{gte}). This involution corresponds to the order 2
rotation around $Q$ in $\R P^n$.

Therefore, we put the projective structure $a'$ on the manifold
$M'$ by gluing together two copies of $(N', c_\la)$ via the
isomorphism $\theta'$ of their boundaries. Then $(M',a')$ admits
an order 2 automorphism $\theta$ which fixes the codimension 2
submanifold $V'$ pointwise and corresponds (under the developing
map) to the order 2 rotation in $\R P^n$.

This concludes Step 1.

\medskip
{\bf Step 2.} Let $N'_\pm$ denote the two copies of $N'$ used to
construct $M'$. Recall that $(N',c)$ is tiled by $m-1$ isometric
copies $O_j$ of the fundamental domain $O$ of $D_m$ (see section
\ref{gte}). The intersection $\cap_j O_j$ is a codimension 2
totally-geodesic submanifold $V\subset N'$.

Let $W_{0}^+=L_0, W_{1}^+,...,W_{m-1}^+\subset N'_+$ denote the
flat hypersurfaces which appear as boundary components of the
domains $O_j\setminus P$. For each $j$ let
$W_{j}^-:=\theta(W_j^+)\subset M'$. Then for each $j$, the flat
hypersurfaces with boundary $W_{j}^-, W_{j}^+$ match in $M'$ to
form a flat hypersurface (without boundary) $S_j\subset M'$ which
is invariant under the automorphism $\theta: (M',a')\to (M',a')$.

The flat hypersurfaces $S_j$ cut the projective manifold $(M',a')$
into components, each of which is (projectively) isomorphic to the
convex hyperbolic manifold with corners $O$.

We now pass to the universal cover $X=(\t{M}', \t{a}')$ of $(M',
a')$. The codimension 2 submanifold $V'$ lifts to $X$ to a
disjoint union of codimension two submanifolds, each of which is
isomorphic to the open $n-2$-disk.

The inverse image $\t{S}\subset X$ of $S:= \cup_j S_j$  is a union
of flat hypersurfaces in $X$ (called {\em walls}) which intersect
along codimension 2 submanifolds above. The closure $C'_j$  of
each component of $X\setminus \t{S}$ is a convex subset, which is
projectively isomorphic to the universal cover of the hyperbolic
manifold with corners $O$. Thus we obtain a covering of $X$ by
closed subsets which are:

1. Codimension 0 strata ({\em cells}): Convex sets $C'_j$.

2. Codimension 1 strata: Facets $F_i'\subset C'_j$, which are
$n-1$-dimensional intersections of $C_j'$ with walls.

3. Codimension 2 strata: Components of the preimage of $P\subset
M'$.

\medskip
The nerve of this decomposition of $X$ is a 2-dimensional cell
complex $Z$, where every 2-cell has $2(m-1)\ge 14$ edges.

By convexity, the links of vertices of $Z$ do not contain any
bigons. Thus $Z$ satisfies the $C'(1/14)$ small cancellation
condition.

\medskip
We now analyze the unions of the adjacent cells. Suppose that
$C_1', C_2'$ share a facet $F'$.

\begin{lem}
$C_1'\cup C_2'$ is isomorphic to a proper convex subset of $\R
P^n$.
\end{lem}
\proof Let $dev: X\to \R P^n$ denote the developing map of $(M',
a')$. Then, by the construction of the projective structure $a'$,
we can assume (after post-composing $dev$ with a projective
transformation) that:
$$
B_i=dev(C_i'), i=1, 2,
$$
are relatively compact convex subsets of $\R^n$ which are
separated by a hyperplane $H\subset \R^n$ passing through the
origin. Let $P$ denote the projective closure of $H$. The
intersection $B_1\cap B_2$ is a facet $\Phi$ contained in $H$. Let
$H_i^+$ denote the closed half-space in $\R^n$ bounded
 by $H$ and containing $B_i, i=1,2$. Let $p\in \R P^n\setminus
 \R^n$be the point dual to $H$. Then there exists $t\in \R$ such
 that the union
 $$
U_t=B_1\cup A_{P,p,t} (B_2)
 $$
is a convex set contained in the unit ball $B$ with center at the
origin. The Euclidean reflection $\tau$ in the hyperplane $H$
preserves the union $U_t$.

Suppose that there exists a point $x_1\in B_1$ such that the
projective line $l$ through $x_1, p$ crosses $H$ in a point $y$
which is not in $\Phi$. Then, by symmetry, $l$ contains a point
$x_2=\tau(x_1)\in A_{P,p,t} (B_2)$. Thus the segment $\ol{x_1
x_2}$ contains a point $y\notin U_t$. Contradiction.

Therefore we can apply Corollary \ref{C1} to $B_1\cup B_2$ (with
$D:=B_1, E:=A_{P,p,t} (B_2)$) and conclude that $B_1\cup B_2$ is
convex. \qed

\medskip
We now {\em de-projectivize} the projective manifold $X$: We
replace each cell $C_j'$ with a convex cone $\t{C}_j'$, etc. The
result is an affine $n+1$-manifold $\t{X}$ which is obtained by
gluing convex cones $\t{C}_j'$. All the conditions of Theorem
\ref{convex} are satisfied by $\t{X}$. It follows from Theorem
\ref{convex} that $\t{X}$ is isomorphic to a proper homogeneous
open convex cone in $\R^{n+1}$. Thus $X$ is isomorphic to a proper
open convex subset of $\R P^n$ and therefore  the manifold $M'$ is
a convex projective manifold. \qed

\begin{cor}
The projective manifold $M'$ is strictly convex.
\end{cor}
\proof Since $M'$ admits a metric of negative curvature,
its fundamental groups are Gromov-hyperbolic. Since $M'$
is convex, it is strictly convex by Theorem
\ref{benoist}. \qed

\bibliography{$HOME/BIB/lit}
\bibliographystyle{siam}

\medskip

Department of Mathematics,

University of California,

Davis, CA 95616, USA,

kapovich$@$math.ucdavis.edu

\end{document}